\documentclass[a4paper,10pt]{amsart}
\usepackage{smooth-pl}
\hypersetup{pdftitle={Smoothing 3D polyhedral spaces}%
,pdfauthor={Nina Lebedeva, Vladimir Matveev, Anton Petrunin and  Vsevolod Shevchishin}
}

\begin{document}

\title{Smoothing 3-dimensional polyhedral spaces}

\author[Lebedeva]{Nina Lebedeva}
\address{N. Lebedeva\newline\vskip-4mm
Steklov Institute, St. Petersburg, Russia.
\newline\vskip-4mm
Mathematical Department of
St. Petersburg State University, Russia.}
\email{lebed@pdmi.ras.ru}

\author[Matveev]{Vladimir Matveev}
\address{V. Matveev\newline\vskip-4mm
Institut f\"ur Mathematik, Friedrich-Schiller-Universit\"at Jena,  Germany.
}
\email{vladimir.matveev@uni-jena.de}

\author[Petrunin]{\\Anton Petrunin}
\address{A. Petrunin\newline\vskip-4mm
Mathematical Department, Pennsylvania State University, USA}
\email{petrunin@math.psu.edu}

\author[Shevchishin]{Vsevolod Shevchishin}
\address{V. Shevchishin\newline\vskip-4mm
National Research University, Higher School of Economics, Moscow, Russia.
}
\email{shevchishin@gmail.com}

\thanks{N.~Lebedeva was partially supported by RFBR grant 
14-01-00062.}
\thanks{A.~Petrunin was partially supported by NSF grant DMS 1309340.}
\date{}
\begin{abstract}
We show that 3-dimensional polyhedral manifolds 
with nonnegative curvature in the sense of Alexandrov
can be approximated by nonnegatively curved 3-dimensional Riemannian manifolds.
\end{abstract}
\maketitle

\section{Introduction}

We define \emph{a polyhedral space} as a complete metric space which admits a locally finite triangulation 
such that each simplex is isometric to a simplex in a Euclidean space.
If in addition the space is homeomorphic to a manifold,
we call it \emph{a polyhedral manifold}.

\begin{thm}{Main Theorem}\label{thm:main}
Assume $P$ is a compact 3-dimensional polyhedral manifold with nonnegative sectional curvature.
Then there is a 
Ricci flow $L^t$ on a 3-dimensional manifold
defined in time interval $(0,T)$
such that $L^t\to P$ as $t\to0$ in the sense of Gromov--Hausdorff
and sectional curvature of $L^t$ is nonnegative for any $t$.

In particular $P$ is a Gromov--Hausdorff limit 
of 3-dimensional nonnegatively curved Riemannian manifolds. 
\end{thm}

\parit{Sketch.}
We introduce so called \emph{$K[\eps]$-pinching}.
A 3-dimensional Riemannian manifold $M$ satisfies $K[\eps]$-pinching
if at any point $x\in M$ and any sectional direction $\sigma_x$ at $x$ we have
\[\Sec(\sigma_x)+\tfrac\eps2\cdot\R(x)\ge 0;\]
here $\Sec(\sigma_x)$ denotes sectional curvature in the sectional direction $\sigma_x$ and $\R(x)$ denotes the scalar curvature at $x$.

In Proposition~\ref{prop:smooth},
we construct a sequence of 3-dimensional Riemannian manifolds $M_n$ converging to $P$ such that
 $M_n$ is  $K[\tfrac1n]$-pinched for each $n$.
 
Further, 
we consider the Ricci flows $M_n^t$ with the initial data $M_n^0=M_n$.
Passing to a limit $L^t$ of $M_n^t$ as $n\to\infty$ we obtain a Ricci flow which is $K[0]$-pinched; i.e., has nonnegative sectional curvature. \qeds

After this work was done, we learned that a similar technique was used by Spindeler in \cite{spindeler} to smooth 3-dimensional Alexandrov spaces of certain type.

\parbf{Acknowledgments.}
We would like to thank 
Yuri Burago,
Bruce Kleiner,
Thomas Richard,
Miles Simon
and Burkhard Wilking.
Big part of the paper was written during our stay at Istanbul Center for Mathematical Sciences,
we want to thank this institute for exelent working conditons.

\section{Remarks and motivations}

\parit{Main motivation.}
This paper was motivated by the  the following conjecture. 

\begin{thm}{Smoothing Conjecture}\label{conj:main}
Assume $P$ is a compact polyhedral space with nonnegative curvature in the sense of Alexandrov 
then it can be presented as a limit of Riemannian  orbifolds\footnote{More formally \emph{as a limit of underlying metric spaces of Riemannian orbifolds}.}
with \emph{nonnegative cosectional curvature}.
\end{thm}

We are about to explain the precise meaning of words \emph{nonnegative cosectional curvature}.
Before that let us state few things.
\begin{enumerate}[(i)]
\item In 3-dimensional case 
\emph{nonnegative cosectional curvature} 
has the same meaning as 
\emph{nonnegative sectional curvature}.
\item In 4-dimensional case 
\emph{nonnegative cosectional curvature} 
has the same meaning as 
\emph{nonnegative curvature operator}.
\item In the dimension $5$ and higher 
any metric with
\emph{nonnegative cosectional curvature} 
also has 
\emph{nonnegative curvature operator},
but converse does not hold.
\end{enumerate}

The curvature tensors on tangent space $T$
form a subspace of $\mathrm{S}^2(\Lambda^2(T))$ which will be denoted by $\mathrm{A}^4(T)$;
here $\mathrm{S}^2(\Lambda^2(T))$ denotes the symmetric square of the space of bivectors of $T$.
Tensor $R\in \mathrm{S}^2(\Lambda^2(T)$ is a curvature tensor 
if and only if it can be presented as a finite sum
\[R=\sum_i\lambda_i\cdot(x_i\wedge y_i)^2,\eqno({*})\]
where $x_i, y_i\in T$ and $\lambda_i\in\RR$;
the latter is equivalent to the 3-cyclic sum identity
on $R$. 

We say that a curvature tensor $R$ has \emph{nonnegative cosectional curvature} if we can find a presentation $({*})$ 
such that $\lambda_i\ge 0$ for all $i$.
We say that $R$ has \emph{positive cosectional curvature} if it can be presented as a sum of the curvature tensor of a round sphere and a curvature tensor  with \emph{nonnegative cosectional curvature}.

Note that the action of $\mathrm{GL}(T)$ extends to the action on $\mathrm{A}^4(T)$.
It turns out that the set of \emph{nonnegative cosectional curvature} forms the minimal closed convex $\mathrm{GL}(T)$-invariant cone which contains the curvature tensor of unit sphere.
By surprising coincidence, the largest proper cone with this property is formed by curvature tensors with nonnegative sectional curvature.

Note that the scalar product on $T$ extends to the scalar product on $\mathrm{A}^4(T)$.
It turns out that the curvature tensor $R\in \mathrm{A}^4(T)$ has nonnegative cosectional curvature if and only if
\[\langle R,S\rangle\ge 0\]
for any tensor $S\in \mathrm{A}^4(T)$
with nonnegative sectional curvature.
This property justifies the term \emph{cosectional}.

The Smoothing Conjecture (\ref{conj:main}) was motivated by the following theorem proved  in  \cite{petrunin}.

\begin{thm}{Theorem}
If a Riemannian manifold $M$ admits a Lipschitz approximation by  polyhedral spaces
with nonnegative curvature in the sense of Alexandrov,
then $M$ has nonnegative cosectional curvature.

Conversely, 
if $M$ is a compact $m$-dimensional Riemannian manifold with positive cosectional curvature
then it admits a Lipschitz approximation by $m$-dimensional polyhedral spaces
with nonnegative curvature.
\end{thm}

Note that Theorem~\ref{thm:main} proves Smoothing Conjecture~\ref{conj:main} for 3-dimensional polyhedral manifolds. 

In fact the proof of Theorem~\ref{thm:main} can be modified to give complete proof of Smoothing Conjecture~\ref{conj:main} for 3-dimensional case.
To do this note that any $3$-dimensional polyhedral $P$
is isometric as a quotient $\bar P/\iota$, where $\bar P$ is a polyhedral manifold and $\iota\:P\to P$ is an isometric involution.
It remains to check that all the constructions in the proof can be made invariant with respect to a given isometric involution.

If a polyhedral space with nonnegative curvature 
admits a smoothing then the link of each simplex has to be homeomorphic to a sphere.
The later follows from Theorem 1.3 proved by Kapovitch in \cite{kapovitch}.
It is expected that this is also a sufficient condition.
Note that not any polyhedral manifold has this property.
For example, let $K$ be the cone over spherical suspension over Poincar\'e homology sphere which can be also thought as
quotient of $\RR^5=\RR^4\times\RR$ by
binary icosahedral group acting on $\RR^4$-factor.
The space $K$ is a topological manifold, but it has edges  with link homeomorphic to Poincar\'e homology sphere.
In particular, $K$ can not be approximated by smooth manifolds with lower curvature bound.

In \cite{boehm-wilking}, B\"ohm and Wilking noted that
nonnegative cosectional curvature survives under Ricci flow.
It gives a hope that smoothing of a polyhedral space 
using Ricci flow used in our proof
can also work in higher dimensions.
Absence of analog of Simon's theorem (see Corollary~\ref{cor:simon}) seems to be the main obstacle on this way.
In particular we do not know the answer to the following question.

\begin{thm}{Question} 
Assume $M$ is a compact $m$-dimensional Riemannian manifold,
$\diam M=D$, $\vol M=v_0$ and
the curvature operator (or cosectional curvature) of $M$ is at least $\kappa$.
Consider the Ricci flow $M^t$ with the initial data $M^0=M$
defined in the maximal interval $[0,T)$.

Is there a positive lower bound for $T$ in terms of $m$, $D$, $\kappa$ and $v_0$? 
\end{thm}

\parit{Smoothing 3-dimensional Alexandrov spaces.}
The following problem goes back to the end of '80-s.
Our paper gives a partial answer.
An other partial answer is given by Spindeler in \cite{spindeler}.
Yet more general problem of that type was considered by Richard in \cite{richard}.

\begin{thm}{Open problem}
Prove that any compact $3$-dimensional Alexandrov space which is homeomorphic to a manifold admits approximation by $3$-dimensional Riemannian manifolds with the same lower curvature bound.
\end{thm}

\section{Preliminaries and notations}

\parbf{Alexandrov's embedding theorem.}
Recall that Alexandrov embedding theorem states in particular that 
any Riemannian metric with curvature $\ge 1$ on $2$-sphere
is isometric to a convex surface in the unit $3$-sphere.
Applying the cone construction 
to the source 
and target spaces 
of this embedding we obtain the following corollary.

\begin{thm}{Corollary}\label{cor:alex}
Let $K$ be a Euclidean cone with nonnegative curvature in the sense of Alexandrov which is homeomorphic to $\RR^3$.
Then $K$ is isometric to the surface of a convex cone in $\RR^4$.

Moreover, if $K$ is smooth away from its tip 
then the surface is smooth away from the tip.
\end{thm}

\parbf{Hamilton's convergence.}
The following statement follows from the main theorem in \cite{hamilton-compactness} and 
the estimate on injectivity radius in terms of diameter, volume, dimension and upper curvature bound.

\begin{thm}{Proposition}\label{prop:ricc-convergence}
Let $M_n^t=(M,g_n^t)$ be a sequence of $m$-dimensional Ricci flows on a compact manifold $M$
defined in the fixed interval $t\in(0,T)$. 
Assume the following two conditions:
\begin{enumerate}[(a)]
\item for each compact interval $\II \subset (0,T)$, 
the curvature of $M_n^t$ is uniformly bounded for all $t\in \II$;
\item there are real numbers $v_0>0$ and $D$ such that 
\[\vol M_n^t\ge v_0\ \ \text{and}\ \ \diam M_n^t\le D\] 
for any $t$ and $n$.
\end{enumerate}
Then, after passing to a subsequence, the solutions converge smoothly to a complete Ricci flow solution $M_\infty^t$, defined for all $t\in (0,T)$.
\end{thm}

\parbf{Simon's theorem.}
The following corollary follows from Theorem 1.9, proved by Simon in \cite{simon-non-collapsed}.

\begin{thm}{Corollary}\label{cor:simon}
Given positive reals $v_0$, $D$ and $\kappa$
there are positive real constants $K$ and $T_0$
such that the following condition holds.

Suppose $M^0=(M,g^0)$ is a compact $3$-dimensional manifold 
such that 
\begin{align*}
\Sec_{M^0}&\ge -\kappa,
&
\vol M^0&\ge v_0,
&
\diam M^0&\le D.
\end{align*}
Let $M^t=(M,g^t)$ be the solution of Ricci flow with initial data $M^0$.
Then $M^t$ is defined in $[0,T_0)$ and 
\begin{align*}
\Sec_M&\ge -K,
&
\vol M^t&>\tfrac{v_0}2,
&
|\Rm_{g^t}|\le \frac Kt
\end{align*}
for any $t\in [0,T_0)$.
Moreover for any two points $x,y\in M$, we have
$$\bigl||x-y|_{g^t}-|x-y|_{g^s}\bigr|
\le 
K\cdot \sqrt{|s-t|},$$
where $|x-y|_{g^t}$ denotes the distance from $x$ to $y$ 
induced by the metric tensor $g^t$. 
\end{thm}

\parbf{Chen--Xu--Zhang pinching.}
Let $(M,g)$ be a compact 3-dimensional Riemannian manifold.
Fix $\eps\ge0$.
We say that $g$ is \emph{$K[\eps]$-pinched} if 
\[\Sec(\sigma_x)+\tfrac\eps2\cdot\R(x)\ge 0\]
for any tangent sectional direction $\sigma_x$ at any $x\in M$.
This condition defines a convex $O(3)$-invariant cone in the space of curvature tensors $\mathrm{A}^4(\RR^3)$.

\begin{thm}{Lemma}\label{lem:pinching}
Let $\eps\ge 0$ and 
$M^t$ be a solution of Ricci flow defined in the interval $[0,T)$.
Assume $M^0$ has $K[\eps]$-pinched curvature
then so is $M^t$ for any $t\in[0,T)$.
\end{thm}

The lemma above is a partial case of the main theorem of Chen, Xu and Zhang in \cite{chen-xu-zhang}.
Lemma 5.1 proved by Simon in \cite{simon} is slightly weaker
but can be used in our proof the same way.

\section{The proofs}\label{sec:rough-smooth}

The proof of Main Theorem will be given in the very end of this section;
it is based on the following proposition.
By $\dist_\GH$ we will denote the Gromov--Hausdorff metric.

\begin{thm}{Proposition}\label{prop:smooth}
Assume $P$ is a compact 3-dimensional polyhedral manifold.
Then there is a real value $\kappa$ and a sequence of 3-dimensional Riemannian manifolds $M_1, M_2,\dots$ 
such that  $M_n$ is $K[\tfrac1n]$-pinched,  $\dist_\GH(M_n,P)<\tfrac1n$ and $\Sec_{M_n}\ge \kappa$ for each $n$.
\end{thm}

Before coming to the proof, we need to discuss structure of singularities of $3$-dimensional polyhedral manifolds.

Let $P$ be as in the proposition.
Assume $x\in P$ is a \emph{singular point};
i.e., $x$ does not have a neighborhood which is isometric to an open subset of $\RR^3$.

The point $x$  will be called  \emph{essential vertex} if the cone at $x$ 
is not isometric to the product $K\times \RR$ where $K$ is a two-dimensional cone.
Note that the set of essential vertices consists of isolated points in $P$ and is therefore finite.

The remaining singular points form a finite number of connected components.
The closure of each component will be called \emph{essential edge}
and the points in the corresponding connected component
will be called \emph{interior points} of this edge.

Consider an essential edge $E$.
Note that a neighborhood $U$ of any interior point of $E$ is isometric to an open subset $U'$ in $K\times \RR$ for some two-dimensional cone $K$.
Under this isometry $\iota\: U\to U'$,
the points on $E$,
are sent to $o_K\times \RR$ in $K\times \RR$,
where $o_K$ denotes the tip of $K$.

Assume $\tau$ is a triangulation of $P$.
Note that any essential vertex of $P$ is a vertex of $\tau$
and any essential edge of $P$ is a union of some edges of $\tau$.
On the other hand a vertex of $\tau$ may not be essential vertex of $P$,
as well as an edge of $\tau$ may not lie in an essential edge of $P$.

Note that each essential edge is a local geodesic.
The essential edge is called \emph{closed} if the geodesic is periodic;
otherwise it is called \emph{open}. 
In the later case the edge connects two vertices 
or a vertex to itself.

Note that the cone $K$ above can be chosen 
the same for all the interior points on $E$.
Assume $\theta_E$ denotes the total angle around the edge;
that is the total angle of the cone $K$.
Then the value $\omega_E=2\cdot\pi-\theta_E$ will be called \emph{curvature} of  $E$. 

Note that by the definition of essential edge,
its curvature can not vanish;
since $P$ has nonnegative curvature,
the curvature of any essential edge is positive.

\parit{Proof of Proposition~\ref{prop:smooth}.}
The construction of the sequence $M_n$ uses two procedures  
(1) the \emph{edge smoothing} and 
(2) \emph{vertex smoothing}.

First we apply \emph{edge smoothing} to $P$.
It produces 
a sequence of manifolds $M_n'$
with isolated singular points for each vertex of $P$
such that 
\begin{enumerate}[(i)]
\item $\dist_\GH(M_n',P')<\tfrac1{2\cdot n}$ for every $n$;
\item The curvature of $M_n'$ is $K[\tfrac1n]$-pinched at any smooth point;
\item Each singular point in $M_n'$ has a conic neighborhood with nonnegative curvature in the sense of Alexandrov.
\end{enumerate}

\parit{Edge smoothing.} 
Assume $E$ is a closed edge.
Denote by $\ell$ its length.
Note that there is 2-dimensional cone $K$, 
a disk $D\subset K$  
and an isometry $\iota\:D\to D$
such that a tubular neighborhood $U$ of $E$ is locally isometric to the space glued from cylinder $D\times [0,\ell]$
by the map $(p,0)\mapsto(\iota(p),\ell)$.

Let us embed $K$ as the graph \[z=k\cdot\sqrt{x^2+y^2}\]
in $(x,y,z)$-space.
Fix a smooth convex even function $\phi(t)$ such that $\phi(t)=|t|$ if $|t|>1$.
Given $\eps>0$, set $\phi_\eps(t)=\eps\cdot\phi(\tfrac t\eps)$.
Denote by $K'_\eps$ the graph
\[z=k\cdot\phi_\eps(\sqrt{x^2+y^2})\]
with induced length metric.

Assume $\eps$ is sufficiently small.
Then there is a rotationally symmetric disk $D'$ in $K'_\eps$
which is isometric to $D$ near the boundary.
Denote by $\iota'\:D'\to D'$ the isometry which coincides with $\iota$ near the boundary of $D$.

Cut the neighborhood $U$ from $P'$ and glue instead $D'\times [0,\ell]/\sim$,
where $\sim$ is the minimal equivalence relation
such that $(p,0)\sim (\iota'(p),\ell)$ 
for any $p\in D'$.
This way we can smooth all the closed edges.

\smallskip

Now assume $E$ is an open edge in $P$.

Note that there is surface $K$
of rotationally symmetric convex cone in $\RR^3$ such that 
$E$ has a neighborhood $\Omega$ 
which is locally isometric to the intersection of a convex neighborhood  of $\{0\}\times(0,\ell)$ in $\RR^3\times\RR=\RR^4$ with $K\times \RR$.

Fix a concave smooth function $f\:[0,\ell]\to\RR$ 
such that for all sufficiently small $t\ge 0$ we have
$f(t)=t$ 
and $f(\ell-t)=t$.

Fix sufficiently small $\eps>0$.
Consider the hypersurface in $\RR^4$ formed by one parameter family of smoothings $K'_{\eps\cdot f(t)}\times\{t\}$;
it also can be described as a graph in $(x,y,z,t)$-space
\[z
=
k\cdot \phi_{\eps\cdot f(t)}
\left(
\sqrt{x^2+y^2}
\right)\]

Let $k_1\le k_2\le k_3$ be the principle curvatures of obtained surface at given point.
For the straightforward choice of functions $\phi$ and $f$,
we have that 
(1) $k_1\le 0\le k_2\le k_3$,
(2) $k_1\cdot k_3\ge \kappa$ for some fixed negative constant $\kappa$ and any $\eps>0$
and (3) $\sup_{k_1\ne 0}\tfrac{k_1}{k_2}\to 0$ as $\eps\to 0$.
In particular, assuming that $\eps$ is sufficiently small, 
the constructed patch has $K[\tfrac1n]$-pinched curvature 
and the sectional curvature at least $\kappa$; here the constant $\kappa$ is independent of $n$.

Note that after the edge smoothing 
the ends of the edge 
have conic neighborhoods with nonnegative curvature.

Applying these operations to all edges of $P$ 
for sufficiently small $\eps=\eps_n>0$ we get the sequence $(M_n')$. 

\medskip

The vertex smoothing produces a Riemannian manifold $M_n$ for each $M_n'$;
it only changes $M_n'$ in a small neighborhood of each vertex leaving this neighborhood nonnegatively curved.

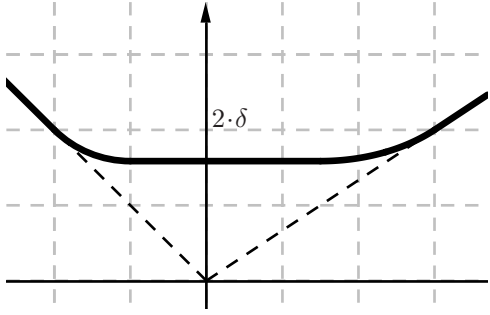
\begin{wrapfigure}{o}{65mm}
\begin{lpic}[t(-0mm),b(-0mm),r(0mm),l(-0mm)]{pics/hat(1)}
\lbl[bl]{27,24;$2{\cdot}\delta$}
\end{lpic}
\end{wrapfigure}

\parit{Vertex smoothing.}
For any singular point,
there is $\eps>0$ such that its $\eps$-neighborhood is conic.
By Corollary~\ref{cor:alex}, this neighborhood
is isometric to an open set in the surface $K$ of convex cone in $\RR^4$.
The surface of cone is smooth at all points except the tip.

We can assume that coordinates $(x,y,z,w)$ in $\RR^4$ are chosen in such a way that $K$ forms a graph $w=f(x,y,z)$ for a nonnegative convex positive homogeneous function $f$.

Fix a convex function $\phi\:\RR_{\ge0}\to \RR_{\ge0}$ which is constant at the points 
$\delta$-close to $0$ 
and identity $2\cdot\delta$-away from zero.  
Note that the graph of composition 
\[\phi\circ f\:\RR^3\to\RR\] 
forms a smooth convex hypersurface in $\RR^4$ which coincides with $K$ sufficiently far from zero.

Cut a conic neighborhood for each vertex of $M_n'$
and glue instead the graph of composition obtained above 
for small enough $\delta>0$. 
This operation produces the needed Riemannian manifold $M_n$.
\qeds

\parit{Proof of Main Theorem.}
Let $M_n$ be the sequence of manifolds provided by Proposition~\ref{prop:smooth}.

Consider the sequence of Ricci flow solutions $M_n^t$ 
with initial data $M^0_n=M_n$.
Applying Corollary~\ref{cor:simon}, 
we get $M_n^t$ are defined in a fixed time interval $[0,T_0)$.

Applying Corollary~\ref{cor:simon} together with Proposition~\ref{prop:ricc-convergence}, 
we can pass to a subsequence of $M_n^t$ which converges to a solution of Ricci flow $L^t$ for  $t>0$.

Since each $M_n$ is $K[\tfrac1n]$-pinched,
Lemma~\ref{lem:pinching} implies that $M_n^t$ is $K[\tfrac1n]$-pinched for any $t$.
It follows that $L^t$ is $K[0]$-pinched;
i.e., $L^t$ has nonnegative sectional curvature for all $t>0$.

By Corollary~\ref{cor:simon} each family $M_n^t$ is uniformly continuous with respect to the Gromov--Hausdorff metric.
Therefore so is the family $L^t$.
In particular $P$ is the Gromov--Hausdorff limit of $L^t$ as $t\to 0$
since it coincides with the limit of $M_n$ as $n\to \infty$.
\qeds

\end{document}